\documentclass[12pt]{article}
\usepackage{amsmath}
\usepackage{amssymb}
\usepackage{dsfont}
\usepackage{cite}
\newsymbol \blackbox 1004
\newcommand{\eh}{\hfill}\newlength{\sperr}

\newenvironment{proof}{{\settowidth{\sperr}{\bf\rm
Proof}%
\par\addvspace{0.3cm}\noindent\parbox[t]{1.3\sperr}
{\bf\rm P\eh r\eh o\eh o\eh f\eh }%
}}{\nopagebreak\mbox{}
$\blackbox$\par\addvspace{0.3cm}}

\def\cls{\mathcal{S}}

\def\nn{\nonumber}
\def\a{\alpha}

\def\G{\Gamma}

\def\Lam{\Lambda}

\def\la{\lambda}
\def\om{\omega}

\def\t{\theta}

\def\vp{\varphi}
\def\vt{\vartheta}
\def\ve{\varepsilon}
\def\wh{\widehat}
\def\wt{\widetilde}
\def\ov{\overline}

\def\p{\partial}

\def\BC{{\mathbb C}}
\def\BR{{\mathbb R}}
\def\BN{{\mathbb N}}

\def\cla{{\mathcal A}}

\def\cll{{\mathcal L}}
\def\clk{{\mathcal K}}

\def\clq{{\mathcal Q}}
\def\clv{{\mathcal V}}

\newcommand{\E}{\mathrm{e}}
\newcommand{\I}{\mathrm{i}}

\newtheorem{Pa}{Paper}[section]
\newtheorem{Tm}[Pa]{{\bf Theorem}}

\newtheorem{Cy}[Pa]{{\bf Corollary}}
\newtheorem{Rk}[Pa]{{\bf Remark}}

\newtheorem{Pn}[Pa]{{\bf Proposition}}

\title{Dynamical and spectral Dirac systems: \\ response function and inverse problems}

\author{A.L. Sakhnovich}

\date{}
\begin{document}
\maketitle

\begin{abstract} We establish simple connections between response functions of 
the dynamical  Dirac systems and $A$-amplitudes and Weyl functions of  the
spectral Dirac systems. Using these connections we propose a new 
and rigorous procedure to recover a general-type dynamical  Dirac system
from its response function as well as a procedure to construct explicit solutions of this
problem.
\end{abstract}



\section{Introduction} \label{intro}
\setcounter{equation}{0} 
Dynamical systems and corresponding control problems are of great interest and three books on dynamical systems \cite{Gros, Has, Kraus} appeared already  in 2015
(see also interesting references therein).  In particular, hyperbolic dynamical systems (to which class the dynamical  Dirac system belongs) are important.
For dynamical systems and related  inverse problems (including inverse problems for dynamical Schr\"odinger and Dirac systems) see, for instance, 
\cite{AMRy, BeMi, Bel2, Byr, GladMo, KuLa, Oks} and numerous references therein.
At the same time  inverse problems for the classical or spectral (frequency domain) self-adjoint Dirac system, also called AKNS, ZS or Dirac type system, had been actively studied since 1950s 
(see \cite{Kre1, LS}). Various interesting results were 
published last years (see, e.g., \cite{AHM, AGKLS, ClGe, EGNT, FKRS3, Hryn, MP, Pu, SaSaR}). 

The interconnections between dynamical  and spectral Dirac systems,
which we establish here, provide a new and rigorous way to recover a general-type dynamical   Dirac system from the response function and to construct explicit solutions
of this problem as well. These interconnections open also further possibilities for the study of both dynamical and spectral Dirac systems.

Spectral Dirac system has the form
\begin{align} &       \label{1.1}
y^{\prime}(x, z )=\I (z j+jV(x))y(x,
z ), \quad
x \geq 0 \quad \left(y^{\prime}(x,z):=\frac{dy}{d x}(x,z)\right),
\end{align} 
where
\begin{align} &   \label{1.2}
j = \left[
\begin{array}{cc}
I_{m_1} & 0 \\ 0 & -I_{m_2}
\end{array}
\right], \hspace{1em} V= \left[\begin{array}{cc}
0&v\\v^{*}&0\end{array}\right],  \quad m_1+m_2=:m.
 \end{align} 
Here $I_{m_i}$ is the $m_i \times
m_i$ identity
matrix and $v(x)$ is an $m_1 \times m_2$ matrix function.

As we already mentioned, dynamical Schr\"odinger and Dirac systems as well as their connections
with response functions and boundary control are of growing interest \cite{AMRa, AMRy, Bel2, BeMi}. Dynamical Dirac system
(Dirac system in 
the {\it time-domain setup}) was studied in the important recent paper \cite{BeMi}.  The dynamical Dirac system considered in \cite{BeMi} is an evolution system of hyperbolic type
and has the following  form:
\begin{align}\label{1.3}
& \I u_t+Ju_x+\clv u=0 \quad (x \geq 0, \quad t \geq 0); \\ & \label{1.3'}
u=\begin{bmatrix}u_1 \\ u_2 \end{bmatrix}, \quad J=\begin{bmatrix} 0 & 1 \\  -1 & 0 \end{bmatrix}, \quad \clv=\begin{bmatrix} p & q \\  q & -p \end{bmatrix}, \quad u_t:=\frac{\p u}{\p t},
\end{align}
where $p=p(x)$ and $q=q(x)$ are real-valued functions of $x$, and initial-boundary conditions are given by the equalities
\begin{align}\label{1.4}
u(x,0)=0, \,\, x \geq 0; \quad u_1(0,t)=f(t), \,\, t \geq 0.
\end{align}
Here $f$ is a complex-valued function (so called boundary control) and the input-output map (response operator) $R: \, u_1(0,\cdot) \rightarrow u_2(0,\cdot)$
is of the convolution form $Rf=\I f +r*f$. The inverse problem consists in recovery of the potential $\clv$ from the {\it response function} $r.$
This  inverse problem was considered in  \cite{BeMi} using boundary control methods.

We note that recent results on dynamical Schr\"odinger and Dirac equations are based on several earlier works.
In his paper \cite{Blag} from 1971, A.S.~Bla-gove\v{s}\v{c}enskii considered dynamical system
\begin{align} \label{h1}&
u_{tt}-u_{xx}+ \clq (x)u_x=0
\end{align}
with boundary control $u(0,t)=f(t)$, and solved inverse problem to recover $\clq$ from $f$. A.S. Blagove\v{s}\v{c}enskii
established important connections between his problem and spectral theory of string equations.
This work was developed further in \cite{Bel1} (see also references therein), where response operator appears
in  inverse problem. Finally, the inverse problem to recover the matrix potential $\clq(x)$ of the dynamical Schr\"odinger equation
\begin{align} \label{h2}&
u_{tt}-u_{xx}+ \clq (x)u=0
\end{align} 
from the response function was considered in \cite{ABI} by  S. Avdonin, M. Belishev, and S. Ivanov. 

Similar to the case of the spectral Dirac and Schr\"odinger equations, the dynamical Dirac  equation is
a more general object than the dynamical Schr\"odinger equation. More precisely, setting
in \eqref{1.3} 
\begin{align} \label{h3}&
p(x)=0, \quad q(x)=g_x(x)/g(x), \quad {\mathrm{where}} \quad g_{xx}(x)=  \clq (x)g(x),
\end{align} 
and rewriting \eqref{1.3}, \eqref{1.3'} in the form
\begin{align} \label{h4}&
(u_1)_t=\I\big((u_2)_x+qu_2\big), \quad (u_2)_t=\I\big(-(u_1)_x+qu_1\big),
\end{align} 
we obtain a dynamical Schr\"odinger equation
\begin{align} \nonumber
(u_1)_{tt}&=(u_1)_{xx}-(q_x+q^2)u_1=(u_1)_{xx}-\big(g_{xx}/g\big)u_1
\\  \label{h5}&
=(u_1)_{xx}-\clq u_1.
\end{align}
For interesting applications of the interconnections between  spectral Dirac and Schr\"odinger equations
see, for instance, the papers \cite{BolGe, GeSS, GeGo, EGNT, EGNST} and references therein.

In our paper, we come from the dynamical to the spectral Dirac system taking Fourier transformation of both
parts of \eqref{1.3}. For that we should be able to estimate the behavior of solutions
$u$ of \eqref{1.3}--\eqref{1.4}. Our next section is dedicated to the necessary estimates on $u$.
A simple connection between response function of the dynamical Dirac system
and Weyl function of the corresponding spectral Dirac system is derived in Section \ref{Resp}.
Recovery of the dynamical Dirac system from the response function is described in Corollary \ref{Cys}
in Section \ref{GTIP}. The connection between response function and so-called $A$-amplitude is described
in Section \ref{GTIP} as well. We note that  interrelations between dynamical and spectral Dirac systems provide 
a procedure for solving  the general-type inverse problem for  dynamical Dirac system
in much greater generality than in \cite{BeMi}.

Explicit methods in spectral theory and in  construction of solutions of equations is an interesting domain, which is actively
developed (to a great extent independently from the general theory). In particular, the methods include
different versions of B\"acklund-Darboux transformations, Crum-Krein transformations and commutation methods
(see, e.g.,  \cite{C1, Ci, Cr, D, Ge, GeT, MS, KoSaTe, KrChr, SaAmmnp10, T00} and various references therein).
In Section \ref{STIP}, we recover explicitly dynamical Dirac system from response function using GBDT (generalized B\"acklund-Darboux transformation) technique from
\cite{GKS6, SaA2, SaAmmnp10, SaSaR} (see also references therein). 

In order to avoid the usage of the same letters for different items, we often use the calligraphic font (e.g., $\clv$)  for functions and operators corresponding
to dynamical system \eqref{1.3}. As usual, $\BR$ stands for the real axis, $\BC$ stands for the complex plane,  $\BC_+$ is the open upper half-plane
$\{z: \,\Im(z)>0\}$ and $\BC_M$ is the half-plane $\{ z: \, \Im (z)>M>0\}$. The notation $\ov{\t}$ stands for the complex conjugate of $\t$
when $\t$ is a complex number, and for the vector with the entries which are complex conjugates of the entries of $\t$ when $\t$ is a vector. 
By  l.i.m. we denote the  entrywise limit of matrix functions in the $L^2$ norm on finite intervals, semiaxes or $\BR$ depending on the context.
We say that the matrix function is boundedly (continuously) differentiable when its derivative is bounded in the matrix norm (continuous).
By $B(H)$ we denote the class of bounded linear operators acting  in the Hilbert space $H$, $I$ is the identity operator, $I_m$ is the $m\times m$ identity matrix 
and $A^*$ stands for the  Hermitian adjoint of an operator (or a matrix) $A$.

\section{Preliminaries and estimates on the solutions $u$ of the dynamical Dirac systems} \label{Prel}
\setcounter{equation}{0} 

According to \cite[Theorem 1]{BeMi}, in the case where
 $p, \, q, \, f$ are continuously differentiable (i.e.,   $p, \, q, \, f \in C^1$) and $f(0)=f^{\prime}(0)=0$, there
is a unique classical solution $u$ of  \eqref{1.3}, \eqref{1.4}  and this solution
admits representation
\begin{align}
\label{2.8}&
 u=u_*^f+w^f; \qquad u_*^f(x,t)=f(t-x)\begin{bmatrix} 1 \\ \I \end{bmatrix},  \quad f(t)=0 \quad (t < 0);\\
\label{2.8'}& w^f(x,t)=w(x,t)= \int_x^t f(t-s)\kappa(x,s)ds  \quad (t\geq x\geq 0), 
\\ \label{2.8''}&
 w(x,t)=0  \quad (x> t\geq 0),
\end{align}
where $\kappa(x,s)$ $(x \leq s)$ is continuously differentiable. In particular, formulas \eqref{2.8} and \eqref{2.8''} yield:
\begin{align} \label{2.f}&
u(x,t)=0 \,\, {\mathrm{for}} \,\, 0 \leq t<x \quad({\mathrm{finiteness \, of \, the \, domain \, of \, influence }}).
\end{align}
Representation \eqref{2.8}--\eqref{2.8''} is proved in \cite{BeMi} using Duhamel formula. 
Moreover, the second term $w^f=w$ in the representation \eqref{2.8} of $u$ admits  (see \cite[p.~6]{BeMi}) an expansion
\begin{align} \label{2.9}&
w(x,t)=\sum_{k\geq 0}\cla^{k+1}u_*^f(x,t) \quad  {\mathrm{for}} \quad t\geq x\geq 0, \\
\label{2.10}& 
\big(\cla g\big)(x,t)=- \cls\big(\clv(x)g(x,t)\big),
\end{align}
 and the operator $\cls $ is given by the right-hand side integrals in   \cite[formulas (1.8)--(1.11)]{BeMi}.
 Namely, for $t \geq x$ we have
 \begin{align} \nonumber
\big(\cls h \big)_1(x,t)=&-\frac{1}{2\sqrt{2}}\left(\int_{\la_1}^{\la_2}(\I h_1+h_2)d \ell -\int_{\la_2}^{\la_3}(\I h_1+h_2)d \ell \right.
\\ & \label{2.10'} \left.
+\int_{\la_1}^{\la_4}(\I h_1-h_2)d \ell           \right), \quad h(x,t)=\begin{bmatrix}h_1(x,t) \\ h_2(x,t)
\end{bmatrix}; \\ \nonumber
\big(\cls h \big)_2(x,t)=&\frac{1}{2\sqrt{2}}\left(\int_{\la_1}^{\la_2}(h_1- \I h_2)d \ell -\int_{\la_2}^{\la_3}(h_1- \I h_2)d \ell \right.
\\ & \label{2.10''} \left.
-\int_{\la_1}^{\la_4}(h_1+ \I h_2)d \ell \right),
\end{align}
 where $\la_i \in \BR^2$ for $i=1,2,3,4$ and $\la_1=(x,t)$, $\la_2=(0,t-x)$, $\la_3=(t-x, 0)$, $\la_4=(x+t,0)$. 
 Here  $\int_{\la_i}^{\la_k} \cdot d\ell$ stands for the integral along the interval $[\la_i, \, \la_k]$ and $\ell$  is the length.
 
 If $h(x,t)=0$ for $t<x$, we have $\big(\cls h\big)(x,t)=0$ for $t<x$,  and therefore,
 in view of \eqref{2.8}, we shall need here
 only the expression for $\cls h$ when $t \geq x$. Moreover, since we assume that $f(t)=0$ for $t<0$,
 it follows from \eqref{2.8''} and \eqref{2.10} that \eqref{2.9} holds for $t<x$ as well.
 
 Let us also assume that $\clv$, $f$ and $f^{\prime}$ are bounded:
\begin{align} \label{2.11}&
\sup_{x >0}\|\clv(x) \|<M_1, \quad \sup_{t>0}\left\|f(t)\begin{bmatrix} 1 \\ \I \end{bmatrix}\right\|<c_0 ,
\quad \sup_{t>0}\left\|f^{\prime}(t)\begin{bmatrix} 1 \\ \I \end{bmatrix}\right\|< \wt c_0.
\end{align}
\begin{Pn} \label{PnEst} Let $p, \, q, \, f$ be continuously differentiable  and let equalities  $f(t)=f^{\prime}(t)=0$ hold for
$t \leq 0$. Assume that \eqref{2.11} is valid. Then the solution $u$ of the dynamical Dirac system \eqref{1.3},
such that \eqref{1.4} and \eqref{2.f} are valid, satisfies the following inequalities
\begin{align} \label{2.!}&
\| u(x,t)\| \leq c_0\E^{Mt}, \quad \| u_t(x,t)\| \leq \wt c_0\E^{Mt}, \quad \| u_x(x,t)\| \leq M_2 \E^{Mt},
\end{align}
where $x\geq 0$ and $t\geq 0$,  $M_2>0$ is some constant, and $M=2 \sqrt{2} M_1$.

Moreover, the functions $u(x,t)$, $u_t(x,t)$ and $u_x(x,t)$ are continuous in the quarterplane $x\geq 0, \, t \geq 0$.
\end{Pn}
\begin{proof}.
Using  \eqref{2.10}-- \eqref{2.10''} and the first inequality in \eqref{2.11}, it is easy to show that 
\begin{align} \nonumber &
\| \big(\cla  g\big)(x,t)\|\leq 2 \sqrt{2} c_k M_1 t^{k+1}/(k+1)
\end{align}
for the case that $\|g(x,t)\|\leq c_k t^k $, $t\geq x$ (recall that we are interested in the case $g(x,t)=0$ for $t<x$). Hence, relation \eqref{2.8} and the second inequality
in  \eqref{2.11} 
imply that 
\begin{align} \label{2.11'}&
\| \big(\cla^{k+1}  u_*^f\big)(x,t)\|\leq c_0 (2 \sqrt{2} M_1)^{k+1} t^{k+1}/(k+1)!.
\end{align}
It follows from \eqref{2.9} and \eqref{2.11'} that
\begin{align} \label{2.12}&
\|w(x,t)\| \leq c_0( \E^{Mt}-1) \quad  (M=2 \sqrt{2} M_1).
\end{align}
Formulas \eqref{2.8}, \eqref{2.11} and \eqref{2.12} yield the inequality $\| u(x,t)\| \leq c_0\E^{Mt}$.

It is easy to see that the functions $\big(\cla^{k}  u_*^f\big)(x,t)$ ($k\geq 0$) are continuous at $x=t$ and $\big(\cla^{k}  u_*^f\big)(t,t)=0$.
Thus, using the first and third inequalities in \eqref{2.11} we (similar to the derivation of formula \eqref{2.11'}) derive
\begin{align} \label{2.12'}&
\| \big(\cla^{k+1}  u_*^f\big)_t(x,t)\|\leq\wt  c_0 (2 \sqrt{2} M_1)^{k+1} t^{k+1}/(k+1)!,
\end{align}
and the second inequality in \eqref{2.!} follows.
Finally, the third inequality in \eqref{2.!} is immediate from \eqref{1.3}, the first inequality in \eqref{2.11}
and the  first and second inequalities in \eqref{2.!}. 

Recall that $u(x,t)$ is continuous at $x=t$. According to \eqref{2.8} and \eqref{2.8'}, the functions $u_x(x,t)$ and $u_t(x,t)$ are continuous at $x=t$.
It is evident that $u$, $u_x$ and $u_t$ are continuous outside $x=t$ as well.
Thus, $u$, $u_x$ and $u_t$ are continuous in the quarterplane.
\end{proof}
Proposition \ref{PnEst} yields the following corollary.
\begin{Cy}\label{Cy2.2} We can apply to $u$ the transformation$:$ 
\begin{align} \label{2.13}&
\wh u(x,z)=\int_{0}^{\infty}\E^{\I z t} u(x,t)dt, \quad z \in \BC_M=\{ z: \, \Im (z)>M\},
\end{align}
where $\wh u$ stands for the Fourier transformation of $u$ $($and Fourier transformation is taken, for the sake of convenience, for the fixed values  $x,z)$. 
Moreover, the same transformation can be applied to $u_t$ and we have
\begin{align} \label{2.13'}&
\I \int_{0}^{\infty}\E^{\I z t} u_t(x,t)dt=z\wh u(x,z), \quad z \in \BC_M.
\end{align}
\end{Cy}
Using the mean value theorem and the fact that $u_x$ is continuous and satisfies the third inequality in \eqref{2.!}, we obtain our next corollary.
\begin{Cy}\label{Cy2.3} Fourier transformation can be applied to $u_x$ and
\begin{align} \label{2.13''}&
\int_{0}^{\infty}\E^{\I z t} u_x(x,t)dt=\frac{d}{d x}\wh u(x,z)=\wh u^{\prime}(x,z), \quad z \in \BC_M.
\end{align} 
\end{Cy}
Now, applying the Fourier transformation to the dynamical Dirac system \eqref{1.3}, we derive
\begin{align} \label{2.16}&
z\wh u(x,z) +J \wh u^{\prime}(x,z) +\clv(x)\wh u(x,z)=0.
\end{align}
We note that various generalized Fourier transformations are successfully  used (see, e.g., \cite{ArD, CoGI, SaSaR}) for solving inverse problems.

An estimate for $\kappa(x,t)$ is also necessary for our further considerations. According to \cite[p. 7]{BeMi}, the formula
\begin{align} \label{p1}&
\kappa(x,t)=\sum_{k \geq 0}\cla^k\big(\cla u_*^{\delta}\big); \quad \big(\cla u_*^{\delta}\big)(x,t)=- \lim_{\ve \to +0}\cls\left((p+\I q)\delta_{\ve}\begin{bmatrix} 1 \\ \I \end{bmatrix}\right);
\\ & \nonumber \delta_{\ve}(t):=1/\ve \,\, {\mathrm{for}} \,\, 0\leq t\leq \ve, \quad   \delta_{\ve}(t):=0 \,\, {\mathrm{for}} \,\, t >\ve 
\end{align}
is valid. Using \eqref{2.10'}--\eqref{2.11} and the second equality in \eqref{p1}, it is easy to show that
\begin{align} \label{p2}&
\| \big(\cla u_*^{\delta}\big)(x,t)\| \leq 2\sqrt{2} M_1.
\end{align}
From the first equality in \eqref{p1} and from  \eqref{p2}, similar to the proof of Proposition \ref{PnEst}, we obtain the following proposition.
\begin{Pn} \label{PnEst2} Let the conditions of Proposition \ref{PnEst} hold. Then
\begin{align} \label{p3}&
\| \kappa(x,t)\| \leq M\E^{Mt}.
\end{align}
\end{Pn} 

\section{Response and Weyl functions} \label{Resp}
\setcounter{equation}{0} 
Recall that the $m_2 \times m_1$ Weyl matrix function (Weyl function) $\vp(z)$ of the spectral Dirac system \eqref{1.1} with the locally summable potential
$V$ is uniquely determined by the inequality
\begin{align}&      \label{2.1}
\int_0^{\infty}
\begin{bmatrix}
I_{m_1} & \vp(z)^*
\end{bmatrix}
Y(x,z)^*Y(x,z)
\begin{bmatrix}
I_{m_1} \\ \vp(z)
\end{bmatrix}dx< \infty ,\quad z \in \BC_+,
\end{align}
where $Y(x,z)$ is the fundamental $m \times m$ solution of \eqref{1.1} normalized by the condition
$Y(0,z)=I_m$ (see \cite{FKRS3}, \cite[Section 2.2]{SaSaR} and some references therein).
Here $\BC_+$ is the open upper half-plane. Weyl functions $\vp$ introduced by \eqref{2.1} are analytic
and contractive in $\BC_+$.

In view of  \eqref{2.8} and  \eqref{2.13}  we see that
\begin{align} \label{2.17}&
\wh u(x,z) =\wh f(z)\E^{\I z x}\begin{bmatrix} 1 \\ \I \end{bmatrix}+\wh w(x,z),  
\end{align}
where (according to \eqref{2.12})
\begin{align} \label{2.18}&
  \|\wh w(x,z)\|=\left\|\int_x^{\infty}\E^{\I z t}w(x,t)dt\right\|
\leq \frac{c_0}{\Im z -M}\E^{-(\Im z -M)x}.
\end{align}
Hence, for $z\in \BC_M$,  the relation $\wh w(x,z)\in L_2^2(0, \infty)$ holds and implies $\wh u(x,z)\in L_2^2(0, \infty)$.
In other words $\wh u$ is the Weyl solution of the Dirac equation \eqref{2.16}. The transformation of the spectral (i.e., frequency-domain)
Dirac system \eqref{2.16} (and its solutions) into the equivalent form \eqref{1.1} is given by the formula 
\begin{align} \label{2.19}&
y=\clk\wh u, \quad v=\I q - p, \quad \clk= \frac{1}{\sqrt{2}} \begin{bmatrix} \I & 1 \\ - \I & 1 \end{bmatrix},
\end{align}
which means that $\clk\wh u$ is the Weyl solution of \eqref{1.1}, \eqref{1.2}, where $m_1=m_2=1$.

Finally, in view of \eqref{1.4}, \eqref{2.8} and \eqref{2.8'} we obtain
\begin{align}& \label{p4}
u_1(0,t)=f(t), \quad u_2(0,t)=\I f(t)+ \int_0^t r(t-s)f(s)ds,
\end{align}
where $r(t)=\kappa_2(0,t)$ is the response function \cite{BeMi}. Taking into account \eqref{p4} and estimates \eqref{2.11} and \eqref{p3}, we have
\begin{align}& \label{2.20}
\wh u_1(0, z)=\wh f(z), \quad \wh u_2(0, z)=\wh f(z)(\wh r(z)+\I).
\end{align}
Now, using \eqref{2.20}, we express $\clk\wh u$ in terms of the normalized fundamental soluton $Y$:
\begin{align}& \label{2.21}
y(x,z)=\clk\wh u(x,z)=\wh f(z)Y(x,z)\clk\begin{bmatrix}  1 \\ \wh r(z)+\I \end{bmatrix}.
\end{align}
Recall that $y=\clk\wh u$ is the Weyl solution, that is, $y\in L^2_2(0,\infty)$ and that the Weyl
function $\vp(z)$ is uniquely determined (see definition  \eqref{2.1}) by the condition 
$$Y(x,z)
\begin{bmatrix}
1 \\ \vp(z)
\end{bmatrix}\in L^2_2(0,\infty).$$ 
Therefore, formula \eqref{2.21} implies that
\begin{align}& \label{2.22}
\vp(z)=\wh r(z)/(\wh r(z)+2\I).
\end{align}
\begin{Pn}\label{Pn1} The response function $r(t)$ of the dynamical Dirac system \eqref{1.3} is connected
with the Weyl function $\vp(z)$ of the corresponding spectral Dirac system \eqref{1.1} $($where
$m_1=m_2=1)$ via equality \eqref{2.22}.
\end{Pn}
We note that, for the classical  case $m_1=m_2=k$, it is often more convenient to introduce Weyl functions $\vp_H(z)$
(see, e.g., Definition 1.51 \cite{SaSaR}) which belong to Herglotz class instead of being contractive. 
(Recall that Herglotz matrix functions in $\BC_+$ are such matrix functions $\psi(z)$ that $\I(\psi(z)^*-\psi(z))\geq 0$
for all $z\in \BC_+$.)
Namely, 
Weyl matrix function $\vp_H(z)$
is uniquely determined by the inequality:
\begin{align} \label{prlm.se.48}&
\int_0^\infty \left[ \begin{array}{lr} I_k &  \I \varphi_H (z )^*
\end{array} \right]
  K  Y(x, z )^*
 Y(x, z )K ^*
 \left[ \begin{array}{c}
I_k \\ - \I \varphi_H (z ) \end{array} \right] dx < \infty , 
\\ \label{th1}&
{K}:=   \frac{1}{\sqrt{2}}       \left[
\begin{array}{cc} I_k &
-I_{k} \\ I_{k} & I_k
\end{array}
\right],  \quad z \in \BC_+.
\end{align}
A simple connection between $\vp$ and $\vp_H$, which easily follows from \eqref{2.1} and \eqref{prlm.se.48}, implies (for $k=1$) that
\eqref{2.22} may be rewritten in the form
\begin{align}& \label{2.23}
\vp_H(z)=\wh r(z)+\I .
\end{align}
Thus, our procedure \cite{SaA4, ALSJSp, SaSaR} to recover potentials $V$ of Dirac systems \eqref{1.1} can be used to recover
potentials $\clv$ of dynamical Dirac systems from response functions. Moreover, even for the case
of continuously differentiable potentials, equality \eqref{2.23} generates important new results, especially
on {\it explicit} recovery of potentials.

{\bf Open Problem.}  Show that equality  \eqref{2.23} holds for much wider classes of potentials $\clv$ then 
continuously differentiable potentials and is also valid for the non-scalar case $m_1=m_2>1$.
\section{Response function, $A$-amplitude \\ and general-type inverse problem}\label{GTIP}
\setcounter{equation}{0}
\subsection{Inverse problem for the spectral Dirac system}
Since $V$ is bounded, the procedure to recover locally bounded potentials from Weyl functions,
which is given in  \cite[Theorem 5.4]{SaA4}, suffices for our purposes. A quick summary of this procedure
(with some functions and operators multiplied, for convenience, by corresponding  constant scalar factors) 
is presented below (and is close to the summary in \cite[Section 2]{FKS}).

First, introduce a family of convolution operators:
\begin{align}& \label{s1}
S_l=\frac{d}{dx}\int_0^l s(x-t) \, \cdot \, dt, \qquad s(x)=-s(-x)^*,
\end{align}
where the index "$l$" in the notation of the operator $S_l$ indicates the space $L^2_k(0,l)$, in which $S_l$ is acting, and
 the matrix function $s(x)$ is associated with Dirac system. More precisely, $s$ is introduced (for $x > 0$)  
via the Weyl function $\vp_H$:
\begin{align}
& \label{s2}
s(x)^*:=  \frac{d}{dx}\left( \frac{\I}{4 \pi}
\E^{ \eta x}
{\mathrm{l.i.m.}}_{a\to \infty} \int_{- a}^{a}\E^{-\I \xi x}
(\xi +\I \eta)^{-2} \varphi_H (\xi +\I \eta)d \xi \right). 
\end{align}
Here $\eta >0$ and $\varphi_H$ is the Weyl function of a spectral Dirac system \eqref{1.1},
where $m_1=m_2=k$ and $V$ is locally bounded on $[0,\infty)$.
Recall that $\vp_H$ is uniquely determined via the inequality \eqref{prlm.se.48}.
The notation l.i.m.  denotes in \eqref{s2} the  entrywise limit in the norm of  $L^2(0,\, \infty)$.
It was proved in \cite{SaA4} that
$$(\xi +\I \eta)^{-2} \varphi (\xi +\I \eta)\in L^2_{k\times k}(-\infty,\, \infty)$$
for every fixed $\eta>0$,  that  l.i.m. on the right-hand side of \eqref{s2} is differentiable (and does not depend on $\eta >0$),
and so $s(x)$ is well-defined.
Moreover,  $s(x)$ is boundedly differentiable on the intervals $(0, \, l)$, $s(+0)=\frac{1}{2}I_k$, and the operators $S_l$ are
bounded and positive (i.e., $S_l>0$)
for all $l>0$. Thus, we have
\begin{align}\label{s3}&
S_l=I+\int_0^l \om(x-t) \, \cdot \, dt>0, \quad
\om(x):=s^{\prime}(x),
\\ & \label{s4}
\om(x)=\om(-x)^*,  \quad s^{\prime}:=\frac{d}{dx}s,
\end{align}
and $S_l, \, S_l^{-1}\in B\Big(L^2_k(0,l)\Big)$.

Next,  introduce the matrix-functions 
 \begin{align}\label{s5}&
\t_1(x)= [I_k \quad 0] Y(x,0)K^*, \quad  \t_2(x)= [0 \quad I_k] Y(x,0)K^*.
\end{align}
According to  \cite[formula (4.16)]{SaA4} we have
 \begin{align} \label{s6}&
\t_2(x)= \frac{1}{\sqrt{2}} \left( [-I_k \quad I_k]
- \int_0^{2x} \om(t)^* S_{2x}^{-1}[  2s(t) \quad I_k]dt \right),
 \end{align}
where $S_{l}^{-1}$ $(l=2x)$ is applied to $[  2s(t) \quad I_k]$ columnwise.
Recall that $Y$ satisfies \eqref{1.1} and is normalized by $Y(0,z)=I_m$, and that $K$
is given by \eqref{th1}. Hence,
\begin{equation} \label{s7}
K Y(x,0)^* j Y(x,0) K^*\equiv J, \quad Y(x,0) K^* J  K Y(x,0)^* \equiv j, \quad
J:=\left[
\begin{array}{cc}
0& I_k  \\ I_k & 0
\end{array}
\right].
\end{equation}
Using \eqref{1.1}, \eqref{s5} and the second relation in \eqref{s7}, we obtain
 \begin{align} \label{s8}&
 v(x)=\I \t_1^{\prime}(x)J\t_2(x)^*.
 \end{align}
 Finally, from \eqref{1.1}, \eqref{th1}, \eqref{s5} and the second relation in \eqref{s7}
we derive the equalities
 \begin{align} \label{s9}&
 \t_1(0)= \frac{1}{\sqrt{2}}  [I_k \quad I_k], \quad \t_1(x)J\t_2(x)^*\equiv 0, \quad \t_1^{\prime}(x)J\t_1(x)^*\equiv 0,
 \end{align} 
 which uniquely determine $\t_1$, assuming that $\t_2$ is already given.
We can now formulate Theorem 5.4 from \cite{SaA4}.
\begin{Tm} \label{saDip}  Let $\vp_H$ be the Weyl function of a spectral Dirac system \eqref{1.1}, where  
$m_1=m_2=k$ and $V$ is locally bounded, that is,
 \begin{align} \label{s10}&
\sup_{0< x<l}\|V(x)\|<\infty \quad {\mathrm{for \,\, any}} \quad l>0.
 \end{align} 
Then $V$ can be uniquely recovered from $\vp_H$
via the formulas \eqref{1.2}, \eqref{s8} and \eqref{s6}, \eqref{s9}.
Here $\om$ on the right-hand side of \eqref{s6}  is obtained from
\eqref{s2} $($after using $\om=s^{\prime})$ and $S_{2x}$ on the right-hand side of \eqref{s6} is given $($using $\om)$ in \eqref{s3}.
 \end{Tm}
\begin{Rk} For the case $m_1=m_2=k=1$, equalities  \eqref{s9} are equivalent to a much simpler relation
 \begin{align} \label{s9!}&
\t_1\equiv - \ov{\t_2}  j,
 \end{align} 
and we recover $\t_1$ from \eqref{s9!} instead of using \eqref{s9}   $($if $\t_2$ is given$)$.
\end{Rk}

In the seminal paper \cite{GeSi} (see also \cite{Si} as well as some references in \cite{GeSi, Si}), the high energy asymptotics of the  Weyl functions of 
Schr\"odinger operators is expressed in terms of the so called $A$-amplitudes.  High energy asymptotics of the Weyl functions
of self-adjoint (spectral) Dirac systems was studied in \cite{ClGe, SaA4} (see also \cite[formula (32)]{SaAIzV}).
According to \cite[formula (3.25)]{SaA4},  the asymptotic equalities
 \begin{align} \label{s11}&
\vp_H(z)=\I I_k+2\I\int_0^l\E^{\I z x}\big(s^{\prime}(x)\big)^*dx+o\left(z\E^{\I z l}\right), \quad |z| \to \infty \quad (l>0)
 \end{align} 
hold in all the angles $c\Im (z)\geq |\Re(z)|$ in $\BC_+$. Formula \eqref{s11} is an analog of the asymptotic expression for the  Weyl function of 
a Schr\"odinger operator  in terms of the  $A$-amplitude. 
\begin{Rk} Since $2\I \big(s^{\prime}(x)\big)^*$ is an analogue of the $A$-amplitude for the Schr\"odinger operator case,
we call this expression the $A$-amplitude of the spectral Dirac system. We note that, in M.G. Krein's terminology,
$\om=s^{\prime}$ is called the accelerant \cite{Kre2} $($see further explanations in \cite{AGKLS}$)$.
\end{Rk}
\subsection{Response function and $A$-amplitude}
Recall that dynamical Dirac  system \eqref{1.3} generates corresponding spectral Dirac system via formulas
\begin{align}\label{s12}
& \clv=\begin{bmatrix} p & q \\  q & -p \end{bmatrix},  \quad V= \left[\begin{array}{cc}
0&v\\  \ov{v} &0\end{array}\right], \quad v=\I q-p.
\end{align}

\begin{Tm} Let $r(t)$ be the response function of a dynamical Dirac system \eqref{1.3} and let the conditions
of Proposition \ref{PnEst} hold. Then $s$ given by formula \eqref{s2} is well-defined and differentiable, and
we have the equality
\begin{align}\label{s13} &
r(t)=2 \I \ov{s^{\prime}(t)},
\end{align}
that is, the response function $r$ of the dynamical Dirac system coincides with the $A$-amplitude of the corresponding spectral Dirac system.
\end{Tm}
\begin{proof}. We rewrite \eqref{2.23} in the form
\begin{align}\label{s14} &
\vp_H(z)- \I= \int_0^{\infty}\E^{\I z t} r(t)dt, \quad  z \in \BC_M .
\end{align}
In view of the equality $r(t)=\kappa_2(0,t)$ and inequality \eqref{p3}, formula \eqref{s14} yields
\begin{align}\label{s15} &
\vp_H(\xi +\I \eta)- \I= {\mathrm{l.i.m.}}_{a\to \infty}\int_0^{a}\E^{\I (\xi +\I \eta) t} r(t)dt, \quad {\mathrm{where}}  \, \eta > M \, {\mathrm{is}} \, {\mathrm{fixed}} 
\end{align}
and $ \xi \in \BR$, that is, l.i.m. is taken in the norm $L^2(-\infty, \infty)$. Hence, putting
\begin{align}\label{s16} &
r(t)=0 \quad {\mathrm{for}} \quad t<0,
\end{align}
and taking inverse Fourier transformation, we derive
\begin{align}\label{s17} &
\E^{-\eta t}r(t)={\mathrm{l.i.m.}}_{a\to \infty}\frac{1}{2\pi}\int_{-a}^a\E^{-\I \xi t}\big(\vp_H(\xi +\I \eta)- \I\big)d\xi
\end{align}
for each fixed $\eta >M$. Thus, it is immediate that on each semiaxis $(-\infty, l)$ we have
\begin{align}\label{s18} &
r(t)={\mathrm{l.i.m.}}_{a\to \infty}\frac{1}{2\pi}\E^{\eta t}\int_{-a}^a\E^{-\I \xi t}\big(\vp_H(\xi +\I \eta)- \I\big)d\xi .
\end{align}
By standard calculations, using \eqref{s18},  we obtain
\begin{align}\nn 
\int_0^x\int_0^y r(t)dtdy&=-\frac{1}{2\pi}\lim_{a\to \infty}\int_{-a}^a\E^{-\I (\xi+\I\eta) x} \, \frac{\vp_H(\xi +\I \eta)- \I}{(\xi +\I \eta)^2}d\xi +C_1x+C_2
\\ \label{s19} &
\\  \nn &=-\frac{1}{2\pi}\E^{\eta x}
{\mathrm{l.i.m.}}_{a\to \infty}\int_{-a}^a\E^{-\I \xi x} \, \frac{\vp_H(\xi +\I \eta)- \I}{(\xi +\I \eta)^2}d\xi +C_1x+C_2,
\end{align}
where l.i.m. is taken on the intervals $(0,l)$ with respect to $x$.

On the other hand, Residue Theorem implies  that, uniformly for $x\geq 0$, the equality
\begin{align}\label{s20} &
\lim_{a\to \infty}\int_{-a}^a\E^{-\I \xi x}(\xi +\I \eta)^{-2}d\xi =-2\pi x \E^{-\eta x}
\end{align}
holds.  Hence, formula \eqref{s2} can be rewritten in the form
\begin{align}
& \label{s21}
s(x)^*=\frac{1}{2}+  \frac{d}{dx}\left( \frac{\I}{4 \pi}
\E^{ \eta x}
{\mathrm{l.i.m.}}_{a\to \infty} \int_{- a}^{a}\E^{-\I \xi x}
 \frac{\varphi_H (\xi +\I \eta)-\I}{(\xi +\I \eta)^{2}}d \xi \right). 
\end{align}
Finally, \eqref{s13} follows from \eqref{s19} and \eqref{s21}.
\end{proof}
\begin{Cy}\label{Cys} Let $r(t)$ be the response function of a dynamical Dirac system \eqref{1.3} and let the conditions
of Proposition \ref{PnEst} hold. Then the function $v(x)$ is uniquely recovered from $r(t)$ using the formula
\begin{align}
& \label{s22}
s(x)=\frac{1}{2}\left(1+\I \int_0^x \ov{r(t)}dt\right), \quad x>0,
\end{align}
and the procedure given in Theorem \ref{saDip} $($more precisely, via the formulas   \eqref{s8} and \eqref{s3}, \eqref{s6}, \eqref{s9!}$)$.
Next, the potential $\clv(x)$ is uniquely recovered from $v(x)$ using \eqref{s12} or, equivalently, using the first equality in \eqref{s12}
and the formulas 
\begin{align}
& \label{sd1}
p=-\Re(v), \quad q=\Im(v).
\end{align}
\end{Cy}
\subsection{Extensions and conditions on the $A$-amplitude \\ and on response function}
From \eqref{s3} and \eqref{s6} it is clear that in order to recover the values of the potential $V$
(of the spectral Dirac system) on $[0,l]$
we need to know $\om(t)$ or, equivalently, $s(t)$ on the interval $[0,2l]$. In the paper \cite{BeMi},
the dynamical Dirac system is considered on the square $\{(x,t): \,  0\leq x \leq T, \,  0\leq t \leq T\}$.
In a similar way, this  yields the necessity to introduce an "extended" problem (and to extend boundary condition
$f$ and response function $r$ on the interval  $[0,2T]$) in order to recover $V(x)$ on $[0,T]$ from $r$.
In other words, from $r$ given on the interval  $[0,T]$, the potential $V$ is recovered on $[0,T/2]$ only.

It is also interesting to compare the conditions on the Weyl function and $A$-amplitude for the spectral Dirac system
with the conditions on the response function for the dynamical Dirac system.
For the case that $v$ is a square matrix function (or a scalar function),
sufficient conditions for a Herglotz function $\vp_H$ to be a Weyl function can be given in terms
of the spectral function \cite{LS, Mar0}, and  the spectral function is connected with $\vp_H$ via Herglotz
representation. On the other hand,  positive operators $S_l$ are also recovered from the spectral function
(see \cite[Theorem 2.11]{SaSaR}), and in this way the conditions from  \cite{LS, Mar0} are related to the condition of positivity
of $S_l$.
The invertibility of the convolution operators (which up to constant factors coincide with our operators $S_l)$ is required in
 \cite{AGKLS, Kre1} and provides the positivity of $S_l$ as well. Finally, for the case of the $m_1\times m_2$ matrix functions
 $v$, sufficient conditions for $\vp(z)$ to be a Weyl function are  local boundedness of the $A$-amplitude and invertibility (positivity)
 of the corresponding operators $S_l$, see \cite[Theorem 2.54]{SaSaR}. 
 We see that the positivity of the operators $S_l$ acting in
 $L^2(0,l)$ is required from $\vp_H$ to be a Weyl function of  a spectral Dirac system.
 The  condition for $r$ to be a response function, which is given in \cite{BeMi}, is again the positivity
 of a structured operator $\cls$ acting in the vector space $L^2_2(0,2T)$, although certain discrepancies in the definition
of  $\cls$ exist (see \cite[pp. 18 and 23]{BeMi}).

We note that necessary and sufficient conditions for the solvability of the inverse problem, which was studied in \cite{Blag}
for system \eqref{h1}, are again formulated in terms of positivity of certain convolution operators (Krein integral operators).

\section{A special type of inverse problem: \\ explicit recovery of the potential}\label{STIP}
\setcounter{equation}{0} 
Explicit construction of the Weyl function (direct problem) and explicit recovery of the potential $V$ of the spectral Dirac system \eqref{1.1} from the Weyl function 
(inverse problem) were
dealt with in \cite{FKRS2013, GKS1, GKS6,  SaAmmnp10}. We shall use the formulation of these results
(for the case $m_1=m_2=k$) from \cite[Subsection 5.1.1]{SaAmmnp10}.

The class of potentials that we recover is called pseudo-exponential \cite{GKS1}.
Each potential from this class is determined by a fixed number $n\in \BN$, by an $n\times n$ matrix
$A$ and by two $n \times k$ matrices $\vt_1$ and $\vt_2$ such that
\begin{equation}       \label{i1}
 A-A^*=\I \Lam(0)j\Lam(0)^*, \quad \Lam(0):=\begin{bmatrix}\vt_1 & \vt_2 \end{bmatrix}.
\end{equation}
Then $V$ corresponding to the triple of parameter matrices $\{A, \,\vt_1, \,\vt_2\}$ is determined by the second equality in \eqref{1.2}
and by the formula
\begin{align}& \label{i2}
v(x)=-2\I  \vt_1^*\E^{\I xA^*}S(x)^{-1}\E^{\I xA}\vt_2,
\end{align}
where $S(x)$ is given by the equations
\begin{align}& \nonumber
S^{\prime}(x)=\Lam(x)\Lam(x)^*, \quad S(0)=I_n; \quad \Lam^{\prime}(x)=-\I A \Lam(x)j, \quad  \Lam(0)=\begin{bmatrix}\vt_1 & \vt_2 \end{bmatrix}.
\end{align}
The formula defining $S$ and $\Lam$ above is easily rewritten in the form
\begin{align}   &    \label{i3}
S(x)=I_n+\int_0^x \Lam(t)\Lam(t)^*dt>0, \\
&    \label{i4}
 \Lam(x)=\begin{bmatrix}\Lam_1(x) \quad \Lam_2(x)\end{bmatrix}
=\begin{bmatrix}\E^{-\I xA}\vt_1 \quad \E^{\I xA}\vt_2\end{bmatrix}.
\end{align}
Weyl functions $\vp_H$  of Dirac systems with pseudo-exponential potentials are rational functions belonging to Herglotz class
(i.e., $\I\big(\vp_H(z)^*-\vp_H(z)\big)\geq 0$ for $z \in \BC_+$). More precisely, the next statement is valid.
\begin{Pn}\label{PnDP} Let $v$ be a pseudo-exponential potential $($i.e., let $v$ admit representation \eqref{i2}-\eqref{i4}, where \eqref{i1} holds$)$.
Then the Weyl function $\vp_H$ of the spectral Dirac system \eqref{1.1}, \eqref{1.2} with this $v$  $(m_1=m_2=k)$ 
 is given by the equality
\begin{equation}       \label{i5}
\vp_H(z)=\I I_k+2\vt_2^*(z I_n -\a)^{-1}\vt_1, \quad \a:=A-\I \vt_1
\big(\vt_1+\vt_2\big)^*.
\end{equation}
\end{Pn}
Vice versa, each proper rational matrix function $\phi(z)$ such that
\begin{align}\label{i6}&
\lim_{z\to \infty}\phi(z)=I_k; \qquad \I\big(\phi(z)^*-\phi(z)\big)\geq 0 \quad {\mathrm{for}} \quad z\in \BC_+
\end{align}
is the Weyl function of a unique spectral Dirac  system with a pseudo-exponential potential.

In order to deduce which class of response functions corresponds to Weyl functions of the form \eqref{i5} $\,(k=1)$,
we take into account the equality
\begin{align}\label{i7}&
\frac{1}{2\pi \I}\int_{\G}e^{- \I z x}(z I_n - {\cll})^{-1}d z=\exp(-\I x{\cll}),
\end{align}
which holds \cite[p. 557]{ASAK} for   anti-clockwise oriented contours $\G$ and matrices ${\cll}$ such that  the spectrum of
${\cll}$ is situated inside $\G$. In view of \eqref{i5}, we rewrite \eqref{s18} in the form
\begin{align}\label{i8} &
r(t)=(-2\I){\mathrm{l.i.m.}}_{a\to \infty}\frac{\vt_2^*}{2\pi \I}\int_{\G_a}\E^{-\I z t}(z I_n -\a)^{-1}dz \, \vt_1,
\end{align}
where $\G_a$ are anti-clockwise oriented contours consisting of the points 
$$\{z: \, |z-\I \eta|=a, \, \Im(z)<\eta\}\cup \{z: \,  -a\leq z-\I \eta \leq a\},$$
$\eta>0$ is sufficiently large and fixed, $\vt_1$ and $\vt_2$ are column vectors ($\vt_1,\vt_2 \in \BC^n$), and l.i.m. may be considered on any finite interval $(0, l)$.
It is immediate from \eqref{i7} and \eqref{i8} that 
\begin{align}\label{i9} &
r(t)=-2\I{\vt_2^*}\E^{-\I t \a}\vt_1.
\end{align}
We easily check also directly that for $r$ of the form \eqref{i9} we, indeed, have 
\begin{align}\label{i10} &
\wh r(z)=2{\vt_2^*}(z I_n- \a)^{-1}\vt_1 \quad (\Im(z)>\|\a\|),
\end{align}
that is,
\begin{align}\label{i11} &
\vp_H(z)=\wh r(z)+\I =\I+ 2{\vt_2^*}(z I_n- \a)^{-1}\vt_1 \quad (\Im(z)>\|\a\|).
\end{align}
Thus, Proposition \ref{PnDP}, relations \eqref{2.23} and \eqref{i10} and uniqueness in the recovery of $v$ from $\vp_H$
(see Theorem \ref{saDip}) imply the following theorem.
\begin{Tm}\label{TmEDIP} Let $r(t)$ be the response function of a dynamical Dirac system satisfying
conditions of Proposition \ref{PnEst} and assume that $r(t)$ admits representation \eqref{i9}, where the matrix
$A=\a+\I \vt_1
\big(\vt_1+\vt_2\big)^*$ and column vectors $\vt_i\in \BC^n$ $(i=1,2)$ satisfy \eqref{i1}. Then the potential $\clv$ of this dynamical Dirac system
is given by the first equality in \eqref{s12} and formulas \eqref{sd1} and \eqref{i2}-\eqref{i4}.
\end{Tm}
We note that identity \eqref{i1} may be rewritten in an equivalent form
$$\a-\a^*=-\I\big(\vt_1+\vt_2\big)\big(\vt_1+\vt_2\big)^*.$$
\begin{Rk}\label{RkES} The requirement $($in Theorem \ref{TmEDIP}$)$ that $r$ should admit representation
\eqref{i9} could be substituted by some requirements  on the Fourier transform  $\wh r$ of $r$
$($more precisely, by the conditions \eqref{i6} on $\phi(z):=\wh r(z)+\I)$. In that case $\phi$
is the Weyl function $\vp_H$ of some spectral Dirac system with a pseudo-exponential potential,
and so matrices $\a$ and $\vt_i$ are recovered from  
$\vp_H$  following the procedure from
\cite[Theorem 5.4]{SaAmmnp10} $($see also references therein$)$.
\end{Rk}
It would be of interest also to construct explicit solutions $u$ of dynamical Dirac systems similar to the construction
of explicit fundamental solutions of the spectral Dirac systems in \cite{FKRS2013, GKS1, GKS6,  SaAmmnp10}.

\bigskip 
\noindent{\bf Acknowledgments.}
 {This research  was supported by the
Austrian Science Fund (FWF) under Grant  No. P24301.}




\begin{thebibliography}{AGKS}
\bibitem{AHM} 
 S.~Albeverio, R.~Hryniv, and Ya.~Mykytyuk,  \textit{Reconstruction of radial Dirac and Schr\"odinger operators from two spectra}, 
 J. Math. Anal. Appl. \textbf{339}, 45--57 (2008).

\bibitem{AGKLS}
D.~Alpay, I.~Gohberg, M.\,A.~Kaashoek, L.~Lerer, and A.\,L.~Sakhnovich,
\textit{Krein systems and canonical systems on a finite interval: accelerants
with a jump discontinuity at the origin and  continuous potentials},
Integr. Equ.  Oper. Theory \textbf{68}, 115--150 (2010). 

\bibitem{ArD}
 D.\,Z.~Arov and H.~Dym, \textit{
Bitangential Direct and Inverse Problems for Systems of Integral and Differential Equations}, 
Encyclopedia of Mathematics and its Aplications (Cambridge University Press,
Cambridge, 2012).

\bibitem{ABI}
S. Avdonin, M. Belishev, and S. Ivanov,  \textit{Boundary control and an inverse matrix problem for the equation} 
$u_{tt}-u_{xx}+V(x)u=0$, Math. USSR-Sb. \textbf{72}, 287--310 (1992).

\bibitem{AMRa}
 S.~Avdonin, V.~Mikhaylov, and K.~Ramdani,  \textit{
 Reconstructing the potential for the one-dimensional Schr\"odinger equation from boundary measurements}, 
IMA J. Math. Control Inform. \textbf{31}, 137--150 (2014).

 \bibitem{AMRy}
 S.~Avdonin, V.~Mikhaylov, and A.~Rybkin,  \textit{The boundary control approach to the Titchmarsh-Weyl m-function}, I: 
 \textit{The response operator and the A-amplitude}, Comm. Math. Phys. \textbf{275},  791--803 (2007). 

\bibitem{Bel1} 
 M. Belishev,  \textit{Wave bases in multidimensional inverse problems}, Math. USSR-Sb. \textbf{67}, 23--42 (1990). 
 
 \bibitem{Bel2} 
  M.~Belishev,  \textit{Boundary control method in dynamical inverse problems--an introductory course}, in:
  \textit{Dynamical Inverse Problems: Theory and Application}, 85--150, CISM Courses and Lectures, Vol. 529 (Springer, Vienna, 2011).
  
  


\bibitem{BeMi}
 M.~Belishev and V.~Mikhailov,  \textit{Inverse problem for a one-dimensional dynamical Dirac system (BC-method)}, 
Inverse Problems \textbf{30}, 125013 (2014).



\bibitem{Blag}  
A.S. Blagove\v{s}\v{c}enskii,  \textit{The local method of solution of the nonstationary inverse problem for an inhomogeneous string} (Russian),
 Trudy Mat. Inst. Steklov. \textbf{115}, 28--38 (1971).

\bibitem{BolGe}  
D.  Boll\'{e}, F. Gesztesy, H. Grosse, W. Schweiger, and B. Simon,  \textit{Witten index, axial anomaly, and Krein's spectral shift function in supersymmetric quantum mechanics}, 
  J. Math. Phys. \textbf{28}, 1512--1525 (1987).
  
  \bibitem{Byr}
 C.\,I.~Byrnes,  D.\,S.~Gilliam,  C.~Hu, and V.\,I.~Shubov,  
\textit{Zero dynamics boundary control for regulation of the Kuramoto--Sivashinsky equation},
 Math. Comput. Modelling \textbf{52}, 875--891 (2010). 

\bibitem{Calo}   
F.~Calogero,  \textit{A new class of solvable dynamical systems}, J. Math. Phys. \textbf{49}, 052701 (2008).  

\bibitem{C1}
   J.\,L.~Cieslinski,
\textit{An algebraic method to construct the Darboux matrix}, {J. Math. Phys.} \textbf{36}, 5670--5706  (1995).

\bibitem{Ci}
{J.\,L.~Cieslinski},  
\textit{Algebraic construction of the Darboux matrix revisited}, {J. Phys. A} \textbf{42}, 404003  (2009).


\bibitem{ClGe}
 S.~Clark and F.~Gesztesy, \textit{Weyl--Titchmarsh M-function asymptotics, local uniqueness results, trace formulas, and Borg-type theorems for Dirac operators},  
 Trans. Amer. Math. Soc. \textbf{354}, 3475--3534 (2002).
 
 \bibitem{Cr}
M.\,M.~Crum, 
\textit{Associated Sturm--Liouville systems}, {Quart. J. Math., Oxford  II Ser.} \textbf{6}, 121--127 (1955).

 
 \bibitem{CoGI}
  A.~Constantin, V.\,S.~Gerdjikov, and R.\,I.~Ivanov,  \textit{Generalized Fourier transform for the Camassa--Holm hierarchy}, 
	Inverse Problems \textbf{23}, 1565--1597 (2007).

\bibitem{D}
P.\,A.~Deift, 
\textit{Applications of a commutation formula}, { Duke Math. J.} \textbf{45}, 267--310 (1978).
  
 \bibitem{EGNT}
 J.~Eckhardt, F.~Gesztesy, R.~Nichols, and G.~Teschl,  \textit{Supersymmetry and Schr\"odinger-type operators with distributional matrix-valued potentials}, 
 J. Spectr. Theory \textbf{4}, 715--768 (2014). 
 
 \bibitem{EGNST}
J. Eckhardt, F. Gesztesy, R. Nichols, A. Sakhnovich, and G. Teschl,  \textit{Inverse Spectral Problems for Schr\"odinger-Type Operators with Distributional Matrix-Valued Potentials},
Differential Integral Equations \textbf{28},   505--522 (2015).

 \bibitem{FKRS3}
B.~Fritzsche, B.~Kirstein, I.\,Ya.~Roitberg, and A.\,L.~Sakhnovich, 
\textit{Recovery of Dirac system from the  rectangular 
Weyl matrix function}, Inverse Problems \textbf{28}, 015010 (2012).

\bibitem{FKRS2013}
B.~Fritzsche, B.~Kirstein, I.\,Ya.~Roitberg, and A.\,L.~Sakhnovich,
\textit{Weyl theory and explicit solutions of direct and inverse problems for a Dirac system with 
rectangular matrix potential},
Oper. Matrices \textbf{7}, 183--196 (2013).

\bibitem{FKS}
B.~Fritzsche, B.~Kirstein,  and A.\,L.~Sakhnovich, \textit{
Weyl functions of Dirac systems and of their generalizations: integral representation,
inverse problem, and discrete interpolation},
J. Anal. Math. \textbf{116}, 17--51 (2012).

\bibitem{Ge}
F.~Gesztesy, \textit{A complete spectral characterization of the
double commutation method}, {J. Funct. Anal.}  \textbf{117}, 401--446 (1993).


\bibitem{GeGo}
F. Gesztesy, J.A. Goldstein, H. Holden, and G. Teschl, 
\textit{Abstract wave equations and associated Dirac-type operators},
Ann. Mat. Pura Appl. \textbf{191}, 631--676  (2012). 

\bibitem{GeSS}
F. Gesztesy, W. Schweiger, and B. Simon, \textit{Commutation methods applied to the mKdV-equation}, Trans. Amer. Math. Soc. \textbf{324}, 465--525 (1991).



\bibitem{GeSi}
F.~Gesztesy and B.~Simon, \textit{A new approach to inverse spectral theory}, II: 
\textit{General real potentials and the connection to the spectral measure}, 
Ann. of Math. (2)  \textbf{152},  593--643  (2000). 

\bibitem{GeT}
F.~Gesztesy and G.~Teschl,  
\textit{On the double commutation method}, {Proc. Amer. Math. Soc.} \textbf{124}, 1831--1840 (1996).


\bibitem{GladMo} 
G.\,M.\,L.~Gladwell and A.~Morassi (eds), \textit{Dynamical Inverse Problems: Theory and Application},
CISM Courses and Lectures, Vol. 529 (Springer, Vienna, 2011).


\bibitem{GKS1}
I.~Gohberg, M.\,A.~Kaashoek, and A.\,L.~Sakhnovich, \textit{Canonical systems with rational spectral densities: explicit formulas and
applications}, Mathematische Nachr. \textbf{194}, 93--125 (1998).



\bibitem{GKS6}
I.~Gohberg, M.\,A.~Kaashoek, and A.\,L.~Sakhnovich, \textit{Scattering problems for a canonical system with a pseudo-exponential potential}, 
 {Asymptotic Analysis} \textbf{29}, 1--38 (2002).

\bibitem{Gros}
 C.~Gros, \textit{Complex and Adaptive Dynamical Systems,  a Primer},  Fourth edition (Springer, Cham, 2015).

\bibitem{Has} 
 Y.~Hasegawa, \textit{Control Problems of Discrete-Time Dynamical Systems}, Second enlarged edition, Studies in Systems, Decision and Control, Vol. 19 (Springer, Cham, 2015).

\bibitem{Hryn}
  R.\,O.~Hryniv,  \textit{Analyticity and uniform stability in the inverse spectral problem for Dirac operators}, J. Math. Phys. 
	\textbf{52},  063513 (2011).
  
  \bibitem{KoSaTe}
{A.~Kostenko, A.~Sakhnovich and  G.~Teschl},  
\textit{Commutation methods for Schr\"odinger operators with strongly singular potentials},  
 Math. Nachr. \textbf{285}, 392--410 (2012). 

 
 \bibitem{Kraus}  
 U.~Krause,  \textit{Positive Dynamical Systems in Discrete Time. Theory, Models, and Applications}, De Gruyter Studies in Mathematics, Vol. 62 (De Gruyter, Berlin, 2015). 
  
\bibitem{Kre1}
M.\,G.~Krein, \textit{Continuous analogues of propositions on polynomials
orthogonal on the unit circle} (Russian), Dokl. Akad. Nauk SSSR
\textbf{105}, 637--640 (1955).

\bibitem{Kre2}
 M.\,G.~Krein, \textit{On the theory of accelerants and S-matrices of canonical differential systems}, 
 Dokl. Akad. Nauk SSSR (N.S.) \textbf{111}, 1167--1170 (1956).
 
 \bibitem{KrChr}
M.\,G.~Krein. 
\textit{On a continuous analogue of a Christoffel formula from the theory of orthogonal polynomials} (Russian), 
Dokl. Akad. Nauk SSSR \textbf{113}, 970--973 (1957).

\bibitem{KuLa}
Y.~Kurylev and  M.~Lassas, \textit{Inverse problems and index formulae for Dirac operators}, 
Advances in Mathematics \textbf{221}, 170--216 (2009).

\bibitem{LS}
B.\,M.~Levitan  and I.\,S.~Sargsjan,  \textit{Sturm--Liouville and Dirac Operators}, Mathematics and its Applications (Soviet Series), Vol. 59
(Kluwer, Dordrecht,  1990).

\bibitem{Mar0}
V.\,A.~Marchenko, 
\textit{Sturm--Liouville Operators and Applications}, 
Operator Theory Adv.  Appl., Vol.  22 (Birkh\"auser, Basel-Boston-Stuttgart, 1986). 

\bibitem{MS}
V.\,B.~Matveev and M.\,A.~Salle, \textit{Darboux Transformations and Solitons} (Springer, Berlin,  1991).


\bibitem{MP}
 Ya.\,V.~Mykytyuk and D.\,V.~Puyda,  \textit{Inverse spectral problems for Dirac operators on a finite interval}, J. Math. Anal. Appl. 
\textbf{386}, 177--194 (2012). 
 
 \bibitem{Oks}
  L.~Oksanen,  \textit{Solving an inverse obstacle problem for the wave equation by using the boundary control method}, 
	Inverse Problems \textbf{29}, 035004 (2013).
 
 \bibitem{Pu}
  D.\,V.~Puyda,  \textit{On inverse spectral problems for self-adjoint Dirac operators with general boundary conditions}, Methods Funct. Anal. Topology \textbf{19}, 346--363 (2013).
  
\bibitem{SaAIzV}
A.\,L.~Sakhnovich, \textit{Asymptotics of spectral functions of an
$S$-colligation}, {Soviet Math. (Iz. VUZ)} \textbf{32}, 92--105 (1988).  
  
\bibitem{SaA2}
 A.\,L.~Sakhnovich, \textit{Dressing procedure for solutions of nonlinear
equations and the method of operator identities}, {Inverse Problems}  \textbf{10},  699--710 (1994).
  
  
\bibitem{SaA4}
A.\,L.~Sakhnovich, \textit{Dirac type and canonical  systems: spectral
and Weyl--Titchmarsh functions, direct and inverse  problems}, Inverse Problems \textbf{18}, 331--348 (2002).  

\bibitem{SaAmmnp10}
A.\,L.~Sakhnovich, \textit{On the GBDT version of the B\"acklund--Darboux transformation and its 
applications to  linear and nonlinear equations and spectral theory},  
Math.~Model.~Nat.~Phenom.~\textbf{5}, 340--389 (2010).
  
 \bibitem{ALSJSp} 
A.\,L.~Sakhnovich, 
\textit{Inverse problem for Dirac systems with locally square-summable potentials and rectangular Weyl functions},  arXiv:1401.3605, Journal of Spectral Theory to appear. 

\bibitem{ASAK}
A.\,L.~Sakhnovich, A.\,A.~Karelin, 
J.~Seck-Tuoh-Mora,   G.~Perez-Lechuga, and
M.~Gonzalez-Hernandez, 
\textit{On explicit inversion of a subclass of operators with  $D$-difference kernels  and Weyl theory of the corresponding canonical systems}, 
Positivity \textbf{14}, 547--564 (2010). 
 
 \bibitem{SaSaR}
A.\,L.~Sakhnovich,  L.\,A.~Sakhnovich, and I.\,Ya.~Roitberg,   \textit{Inverse Problems and Nonlinear Evolution Equations. 
 Solutions, Darboux Matrices and Weyl--Titchmarsh Functions}, {De Gruyter Studies in Mathematics}, Vol.  {47} (De Gruyter, Berlin, 2013).

\bibitem{Si} 
 B.~Simon, \textit{Schr\"odinger operators in the twentieth century},  J. Math. Phys. \textbf{41}, 3523--3555 (2000).  
 
 \bibitem{T00}
 G.~Teschl, 
 \textit{Deforming the point spectra of one-dimensional Dirac operators},  {Proc. Amer. Math. Soc.}  \textbf{126}, 2873--2881  (1998). 

\end{thebibliography}
\end{document}